\documentclass [hidelinks,12pt,openbib]{article}

\usepackage{cancel}
\usepackage{xcolor}
\usepackage[hidelinks]{hyperref}
\usepackage[normalem]{ulem}
 
\usepackage{mathtools}
 \usepackage{tikz}
 \usepackage{cancel}
\usepackage{fixmath}
\usepackage{dsfont}
\usepackage{caption}
\usepackage{pslatex}
\usepackage{bm}
\usepackage{graphicx}
\usepackage{amsmath,amssymb,amsfonts,amsthm}
\usepackage{xcolor}
\usepackage{mathtools}
\usepackage{dsfont}
\usepackage{cite}
\usepackage{ulem}
\usepackage{fancybox,fancyvrb}
\usepackage{framed} 

\usepackage{MnSymbol,wasysym}
\usepackage[export]{adjustbox}

\newtheorem{proposition}{Proposition}

\usepackage{enumitem}
\usepackage{color}
\usepackage{fancyhdr}
\usepackage{esint}
\usepackage{mathtools}
\usepackage[margin=1in]{geometry}
\usepackage{mdframed}

\def\XXint#1#2#3{{\setbox0=\hbox{$#1{#2#3}{\int}$}
     \vcenter{\hbox{$#2#3$}}\kern-.5\wd0}}

\def\E1{\text{E}_1}

\def\L1loc{{L^1_{\rm loc}}}

\def\ig{\includegraphics}

\def\bc{\begin{center}}
\def\ec{\end{center}}

\def\begfr{\begin{mdframed}[
	backgroundcolor=gray!15!white,
	linecolor=blue,
	linewidth=1pt,
	align=center,
	userdefinedwidth=4.99in]}
\def\begfrwhite{\begin{mdframed}[
	backgroundcolor=white!15!white,
	linecolor=black,
	linewidth=0.5pt,
	align=center,
	userdefinedwidth=4.99in]}
\def\begfrblue{\begin{mdframed}[
	backgroundcolor=blue!15!white,
	linecolor=blue,
	linewidth=1pt,
	align=center,
	userdefinedwidth=4.99in]}

\def\endfr{\end{mdframed}}
\def\endfrred{\end{mdframed}}
\def\begfrred{\begin{mdframed}[
	backgroundcolor=red!15!white,
	linecolor=blue,
	linewidth=1pt,
	align=center,
	userdefinedwidth=4.99in]}
\def\endfrblue{\end{mdframed}}

\def\begeq{\begin{equation}}
\def\endeq{\end{equation}}

\def\bbr{\begin{bmatrix*}[r]}
\def\ebr{\end{bmatrix*}}
\def\bb{\begin{bmatrix*}[r]}
\def\eb{\end{bmatrix*}}
\def\bbc{\begin{bmatrix}}

\def\ebc{\end{bmatrix}}

\def\sign{{\rm sign}}

\definecolor{darkmagenta}{rgb}{0.55, 0.0, 0.55}

\definecolor{darkgreen}{rgb}{0.0,0.6,0.0}

\def\R{\mathbb{R}}

\def\0u{\underline{0}}

\def\begeq{\begin{equation}} 
\def\endeq{\end{equation}}

\def\bcenter{\begin{center}}
\def\ecenter{\end{center}}
\def\beq{\begin{equation}}
\def\eeq{\end{equation}}
\def\bmatr{\begin{bmatrix*}[r]}
\def\ematr{\end{bmatrix*}}
\def\bmatc{\begin{bmatrix}}
\def\ematc{\end{bmatrix}}

\usepackage{amsmath}               
  {
      \theoremstyle{plain}
      
  }
  {
      \theoremstyle{plain}
      
  }
   {
      \theoremstyle{plain}
      
  }
   {
      \theoremstyle{plain}
      
  }
  {
      \theoremstyle{plain}
      
  }
    {
      \theoremstyle{plain}
      
  }
\newcommand{\vertiii}[1]{{\left\vert\kern-0.25ex\left\vert\kern-0.25ex\left\vert #1 
    \right\vert\kern-0.25ex\right\vert\kern-0.25ex\right\vert}}
 \usepackage{amsthm}
\newtheorem*{proposition*}{Proposition}
\newtheorem*{corollary*}{Corollary}
\newtheorem*{definition*}{Definition}
\newtheorem*{example*}{Example}
\newtheorem*{lemma*}{Lemma}
\newtheorem*{theorem*}{Theorem}
\newtheorem*{observation*}{Observation}
\newtheorem*{remark*}{Remark}
\usepackage{sectsty}
\sectionfont{\fontsize{15}{10}\selectfont}

\newtheorem*{exmp*}{Example}
\usepackage{stmaryrd}

\begin{document} 

\bc
{\Large \bf A particle method for  continuous \\ \vskip 5pt
Hegselmann-Krause opinion dynamics}
\ec

\vskip 10pt

\bc
{\large 
Bruce Boghosian$^\ast$, Christoph B\"orgers$^\ast$, Natasa Dragovic$^\ast$,  \\ \vskip 3pt Anna Haensch$^{\ast \dagger}$, and Arkadz Kirshtein$^\ast$
}
\ec

\noindent
$^\ast$ Department of Mathematics, Tufts University, Medford, MA 02155

\noindent
$^\dagger$ Data Intensive Studies Center, Tufts University, Medford, MA 02155

\begin{quote}
{\small
\noindent
{\bf Abstract.} 
We derive a differential-integral  equation akin to  the Hegselmann-Krause model of opinion dynamics [R.\ Hegselmann and 
U.\ Krause, JASSS, vol.\ 5, 2002], and propose a particle method for solving the equation. Numerical experiments demonstrate second-order convergence of the
method in a weak sense. We also show that our differential-integral equation can equivalently be stated as a system of differential 
equations. An integration-by-parts argument that would typically yield an {\em energy dissipation} inequality in physical problems then yields 
a  {\em concentration} inequality, showing that a natural measure of concentration increases monotonically. 
}
\end{quote}

\noindent
{\bf keywords:} opinion dynamics; Hegselmann-Krause model; bounded confidence model; particle method


\section{Introduction.} 
People's opinions and beliefs are influenced in complex ways by families, friends, colleagues, media, as well
as politicians and other mega-influencers \cite{Acemoglu_opinion,Bisin_2000,Bisin_2001,Boyd_1985,Cavalli_1981}.
In recent decades, attempts have been made to understand aspects of this process using mathematical modeling and computational simulation; for surveys on opinion dynamics, see for instance
  \cite{Aydogdu_2017,Anderson_2019,Lorenz_2007,Mossel_2017,Proskurnikov_2017,Proskurnikov_2018}.

Many models of opinion dynamics are based on the assumption  that we are influenced more easily by people
whom we {\em almost} agree with to begin with than by those whose views starkly differ from ours.
A similar but more general phenomenon is known as {\em biased assimilation} among psychologists --- our tendency 
to filter and interpret information in such a way that it supports our preconceived notions \cite{Lord_1979}. Models of opinion dynamics
 based on this assumption are known as {\em bounded confidence models} \cite{Aydogdu_2017,Canuto_2012,Mirtabatabaei_2012}. 
A popular example  is due to Hegselmann and Krause \cite{hegsel_krause_2002,Hegselmann_Krause_2005}, 
building on earlier work by Krause \cite{krause_1997,krause_2000}. It has been studied extensively
in the literature (see
for instance \cite{Lorenz_2005, Lorenz_2006, Blondel_2009}), and will be our starting point
here. The Weissbuch-Deffuant model \cite{Weissbuch_Deffuant} is very close to that of Hegselmann and Krause; while
Hegselmann and Krause assume that each opinion holder responds to all nearby opinions simultaneously, 
Weissbuch and Deffuant assume random encounters between pairs of opinion holders with similar views.  For  other bounded
confidence models, see  
\cite{Ben_Naim_2005, Fortunato_2005, Urbig_2003}. 

The original Hegselmann-Krause model is discrete
in both time and  opinion space. Similar models that are
continuous in time \cite{Piccoli_2021}, 
opinion space \cite{Wedin_2015}, or opinion space and time \cite{Goddard_2022} have been proposed as well. We are particularly
interested in fully continuous models, since we plan, in future work, to 
explore the response of 
candidates to a  dynamic electorate. In a
previous paper, we have already discussed the response of candidates
to a {\em static} electorate \cite{Borgers_et_al_candidate_dynamics}. 
We want to describe candidate dynamics in opinion space by 
ordinary differential equations, and find that easiest to do in clean and natural ways if
the opinion dynamics of the electorate are described fully continuously.

We note that ``continuous" does not mean the same to all authors in this field. For instance, in Lorenz's earlier 
papers \cite{Lorenz_2005, Lorenz_2006}, the dynamics 
are discrete in both opinion space and time. The word ``continuous" appears in the titles of 
both papers, but it indicates merely
that the opinions can take arbitrary real values. To us, by contrast, a ``fully continuous" model is one in which a continuum of agents 
changes opinions in continuous time. We note that fully continuous models in our sense were studied by Lorenz 
in \cite{Lorenz_2007}.

In this paper, we derive a fully continuous version of the Hegselmann-Krause model. We start with a
time-continuous, space-discrete model. 
In contrast with many of the existing time-continuous 
models  \cite{Piccoli_2021}, 
we don't interpret particles as agents, but as agent {\em clusters} of different sizes.
Our time-continuous model has a natural space-time-continuous analogue, a differential-integral equation. The time-continuous model that we start out with
can then be interpreted as a numerical method
for the differential-integral 
equation, a particle method. We numerically test  the speed of
convergence of this method. We note that  particle
methods are a natural choice for the numerical simulation of  bounded confidence models because
biased assimilation tends to result in  the formation of
{\em clusters} of like-minded individuals --- groups of friends confirming and equalizing each others' opinions on Facebook or over dinner,
for instance --- causing accuracy issues for numerical methods based on fixed
grids.

We also observe that the differential-integral equation can be translated into a system of partial differential equations
without any integrals, somewhat reminiscent of the Poisson-Nernst-Planck model of electro-diffusion: The density
(of individuals in opinion space, or of charged particles in physical space) moves in a velocity field that itself is determined by
the density via Poisson-like differential equations. In our model, we show that an integration-by-parts argument that 
would lead to an {\em energy dissipation} inequality in physical systems leads to a {\em mass concentration} inequality here.

\section{A time-continuous model.} 
\label{sec:particle_model}

\subsection{Opinion space and opinion holder distributions.}  We assume that any individual's opinions can be characterized by a single real number $x$. In politics, one could think of 
this as the ``left-right axis", with values of $x$ on the left side of the axis corresponding to ``left" views, and values on the 
right side to ``right" ones. This is a gross simplification that captures some aspect of the truth, since political 
views on different issues are correlated: If you 
tell us your thoughts about immigration policy and about allowing  organized prayer in schools, we cannot be sure how you feel about  a single-payer healthcare system, but we do have a guess.

Tony Blair has suggested a re-interpretation of the one-dimensional axis as ``open" (in favor of immigration, multi-culturalism, globalism) and ``closed" (in favor of restricting immigration, culturally
conservative,  primarily focused on one's own country) \cite{Tony_Blair}. The  interpretation of the one-dimensional axis does not, of course,
affect our abstract modeling. 

We refer to the $x$-axis as ``the opinion axis" or ``opinion space". A space-continuous model typically uses a 
time-dependent {\em density} of opinion holders, 
$$
f(x,t), ~~~ x \in \R, ~~~ t \geq 0.
$$
The time $t$ could still tick discretely in such a model, but we are primarily interested in  space-time-continuous models 
in which $t$ flows continuously. We always assume $f(x,t) \geq 0$ and normalize so that 
$$
\int_{-\infty}^\infty f(x,t) ~\! dx= 1 ~~~\mbox{for all $t$}.
$$
More generally and abstractly, the opinion holder distribution $f$ could be a time-dependent Borel probability 
measure on $\R$; however, the only measures without densities that we'll talk about in this paper are weighted
sums of Dirac measures.

\subsection{Particle representation of opinion distributions.} 
Let $X_1,X_2,\ldots,X_n \in \R$, and assume for now that the $X_i$ are the only opinions represented in the electorate. 
If $w_i$ is the fraction of individuals who hold opinion $X_i$, then the ``density" of opinions altogether is the distribution
\begin{equation}
\label{eq:weighted_sum_of_diracs}
 \sum_{i=1}^n w_i \delta(x-X_i),
\end{equation}
where $\delta$ denotes the Dirac delta distribution. The condition that this be a probability measure becomes
$
\sum_{i=1}^n w_i =1.
$

Any Borel probability measure $\mu$ on the real line can  be approximated arbitrarily well, in the distributional sense,
by a weighted sum of delta functions in the form (\ref{eq:weighted_sum_of_diracs}). In fact, let $m \geq 1$ be an integer,
$\Delta x >0$ a real number, and
define, for all integers $i$ with $-m+1 \leq  i  \leq m-1$, 
\begin{equation}
\label{weights}
w_i = \frac{\mu \left[ \left( i - \frac{1}{2} \right) \Delta x, \left( i + \frac{1}{2} \right) \Delta x \right) }{\displaystyle{\sum_{k=-m+1}^{m-1} } \mu \left[ \left( k - \frac{1}{2} \right) \Delta x, \left( k + \frac{1}{2} \right) \Delta x \right) }.
\end{equation}
Then 
\begin{equation}
\label{dirac_sum}
\sum_{i=-m+1}^{m-1}  w_i   \delta \left( x - i \Delta x \right), 
\end{equation}
converges weakly to $\mu$ if $m \rightarrow \infty$ and $\Delta x \rightarrow 0$ in such a way 
that $m \Delta x \rightarrow \infty$.

Any weighted sum of delta functions in the form (\ref{eq:weighted_sum_of_diracs}) can in turn be approximated 
arbitrarily well by 
a smooth density. For instance, the smooth probability density
\begin{equation}
\label{smooth}
\sum_{i=1}^n w_i  \frac{e^{- (x-X_i)^2/(2 \sigma^2)}}{\sqrt{ 2 \pi \sigma^2}} 
\end{equation}
converges to (\ref{eq:weighted_sum_of_diracs}), in the distributional sense, as $\sigma \rightarrow 0$.

\subsection{The dynamics of conformist opinion holders.} 
We assume that the opinion holders in the $i$-th cluster, that is, opinion holders
with opinion $X_i$, consider a weighted average of the opinions of others, in the form
$$
\frac{\sum_{j=1}^n \eta(|X_i-X_j|) w_j X_j}{\sum_{j=1}^n {\eta(|X_i-X_j|) w_j}}
$$
where 
$$
\eta: ~~  [0,\infty) \rightarrow [0,1]
$$
is a decreasing function with $\lim_{z \rightarrow \infty} \eta(z) = 0$, called the 
{\em interaction function}. The further $X_i$ is removed from $X_j$, the less will the $j$-th cluster affect the 
opinion of the $i$-th cluster. 
We then assume that $X_i$ drifts towards a weighted average of opinions (including their own), where nearby opinions 
are weighed more strongly than ones far from $X_i$: 
\begin{equation}
\label{time_evolution_0}
\frac{dX_i}{dt} = \alpha \left( \frac{\sum_{j=1}^n \eta( \left| X_i-X_j \right|) w_j X_j}{\sum_{\ell=1}^n {\eta( \left| X_i-X_\ell \right|)w_\ell}} - X_i \right).
\end{equation}
Here $\alpha >0$ is a parameter determining how eager the opinion holders are to fall in line with those who already
hold opinions similar to theirs. We will take $\alpha=1$. This is just a matter of choosing time units. Using this, and simplifying
a bit, (\ref{time_evolution_0}) becomes
\begin{equation}
\label{time_evolution} 
\frac{dX_i}{dt} =\sum_{j=1}^n ~\!  \frac{ \eta( \left| X_i-X_j \right|) w_j }{\sum_{\ell=1}^n {\eta( \left| X_i-X_\ell \right|)w_\ell}} ~\! (X_j-X_i).
\end{equation}
The model is most closely analogous to that
of Hegselmann and Krause if the interaction function $\eta$ is taken to be the indicator function of an interval $[0, \epsilon]$ with $\epsilon>0$. However, we use
$$
\eta(z) = e^{-z/ \nu}
$$
where $\nu>0$ is a parameter determining how broad-minded the opinion holders are. Larger $\nu$ means greater broad-mindedness.

Setting $w_j=1$ in (\ref{time_evolution}), our
equation simplifies to \cite[eq.\ (6)]{Aydogdu_2017}. 
A time-discrete version of the model of \cite{Aydogdu_2017} also appears in \cite[eqs.\ (3) and (4)]{Canuto_2012}.
A time-discrete model including
 weights can be found for 
 instance in \cite[eq.\ (2)]{Blondel_2009}.  If $\eta$ is taken to be an indicator function, our model 
 becomes a time-continuous version of that in \cite{Blondel_2009}.

\subsection{Examples.} 
\begin{figure}[h!]
\bc
\ig[scale=0.5]{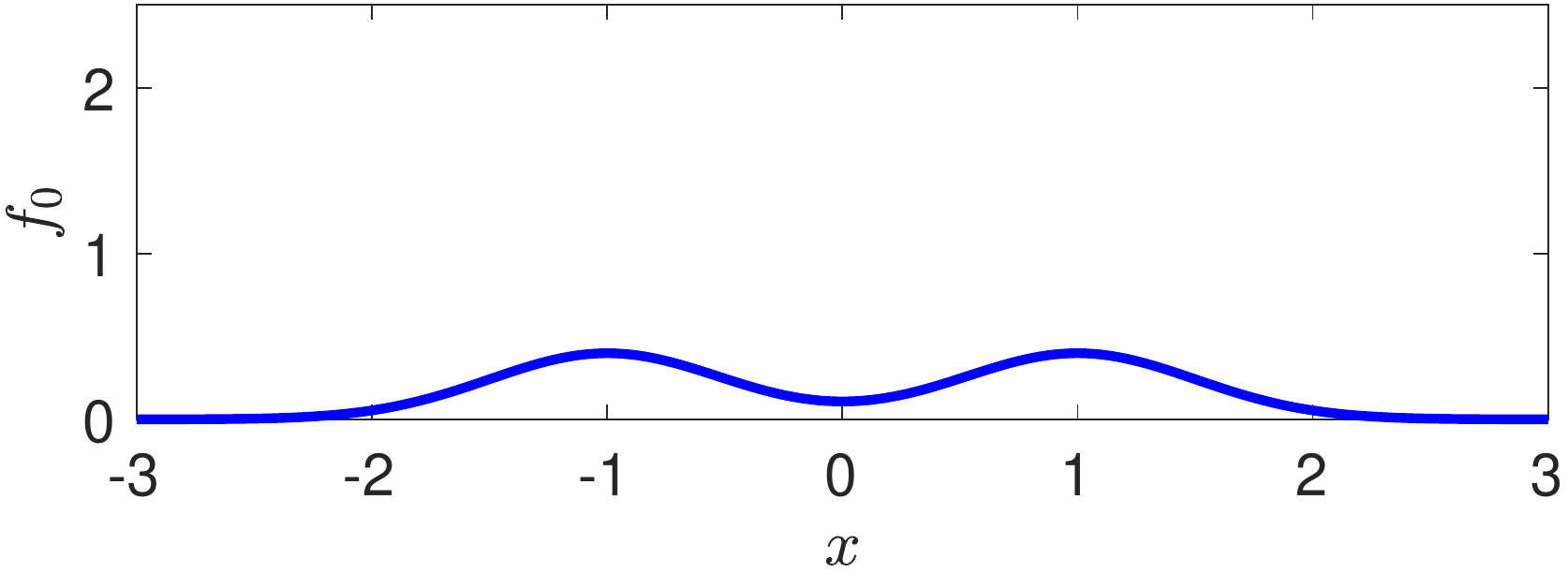}
\caption{An opinion distribution with two distinct ``camps", a ``left" one and a ``right" one.} 
\label{fig:TWO_BUMPS} 
\ec
\end{figure}

We first assume that the initial opinion distribution has the density
$$
f_0(x) = \frac{1}{2} \left( \frac{e^{-2(x+1)^2}}{\sqrt{ \pi/2}} + \frac{e^{-2(x-1)^2}}{\sqrt{ \pi/2}} \right).
$$
The graph of this function is shown in Fig.\ \ref{fig:TWO_BUMPS}. 

\begin{figure}[h!]
\bc
\ig[scale=0.35]{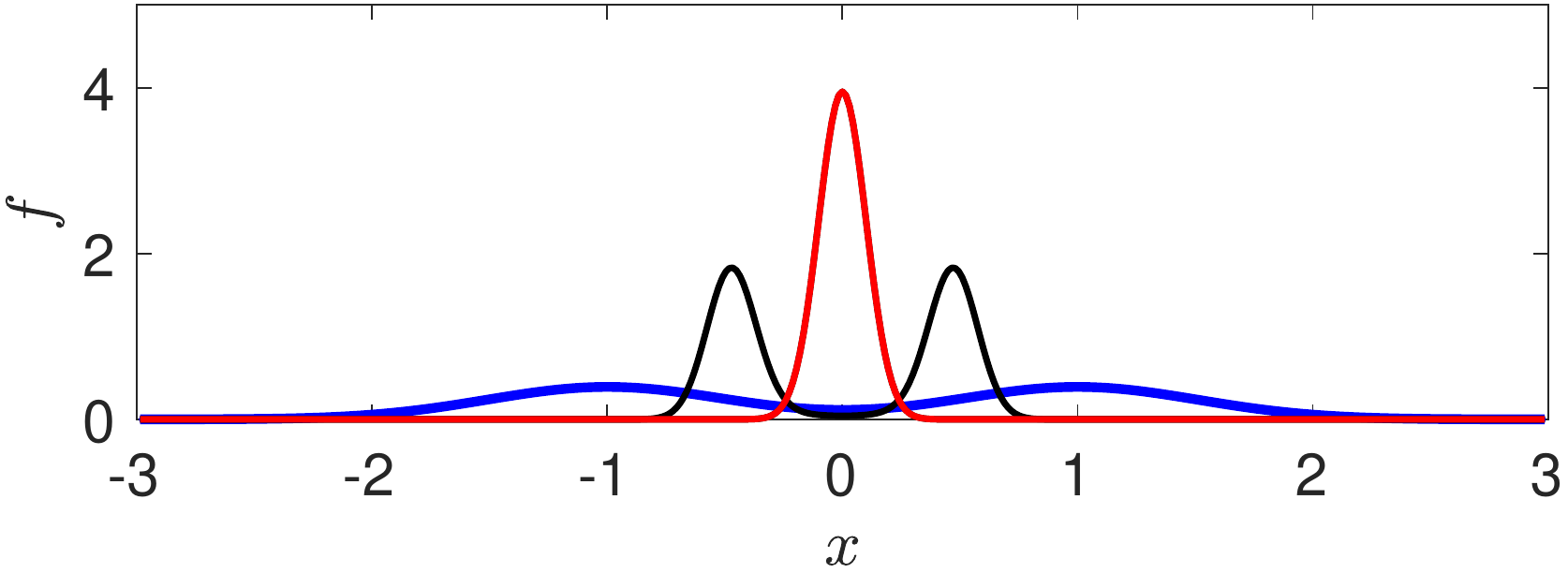}

\bc
\ig[scale=0.35]{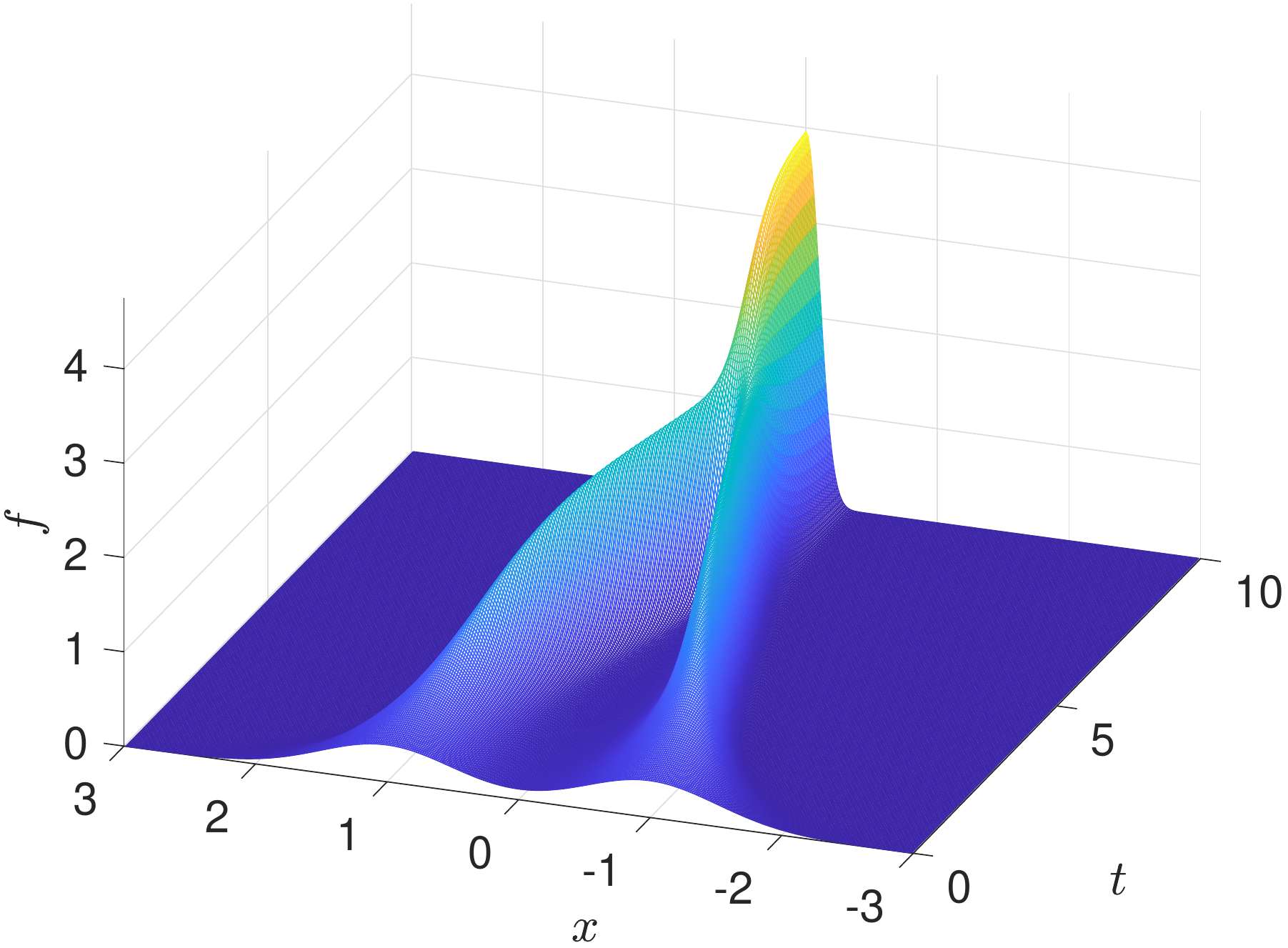}
\ec

\caption{If the initial opinion distribution is that shown in  Fig.\ \ref{fig:TWO_BUMPS}, the initial clusters tighten 
at first, but also move towards one another, and eventually merge. The upper panel shows $f$ at times $0$ (blue), 
5 (black), and 10 (red). The lower panel shows the time evolution as a 3D surface plot.}
\label{fig:TWO_BUMP_EVOLUTION}
\ec
\end{figure}

We approximate this distribution by a weighted sum of 399 Dirac delta functions, 
as described by equations (\ref{weights}) and (\ref{dirac_sum}) with 
$m=200$ and $\Delta x =3/m$, approximating 
$$
\mu \left[ \left( k - \frac{1}{2} \right) \Delta x, \left( k + \frac{1}{2} \right) \Delta x \right)  = \int_{(k-1/2) \Delta x}^{(k+1/2) \Delta x} f_0(x) dx
$$
by 
$$
f_0 \left( k \Delta x \right) \Delta x.
$$
We compute the time evolution  as described by eq.\ (\ref{time_evolution}), 
using $\eta(z) = e^{-2 z}$, and using the midpoint method with $\Delta t = 0.04$.
At each time $t>0$, this results in a weighted sum of Dirac delta functions approximating the opinion distribution. We approximate this sum by a smooth probability density as defined in (\ref{smooth}) 
with $\sigma = 0.1$. The upper panel of Fig.\ \ref{fig:TWO_BUMP_EVOLUTION} shows the resulting densities at times $t=0$ (blue), $5$ (black), and 10 (red). The lower panel shows the same time evolution as a surface plot. The two initial clusters tighten, but they also move towards each other, and eventually they merge into one cluster at the center.

\begin{figure}[h!]
\bc
\ig[scale=0.35]{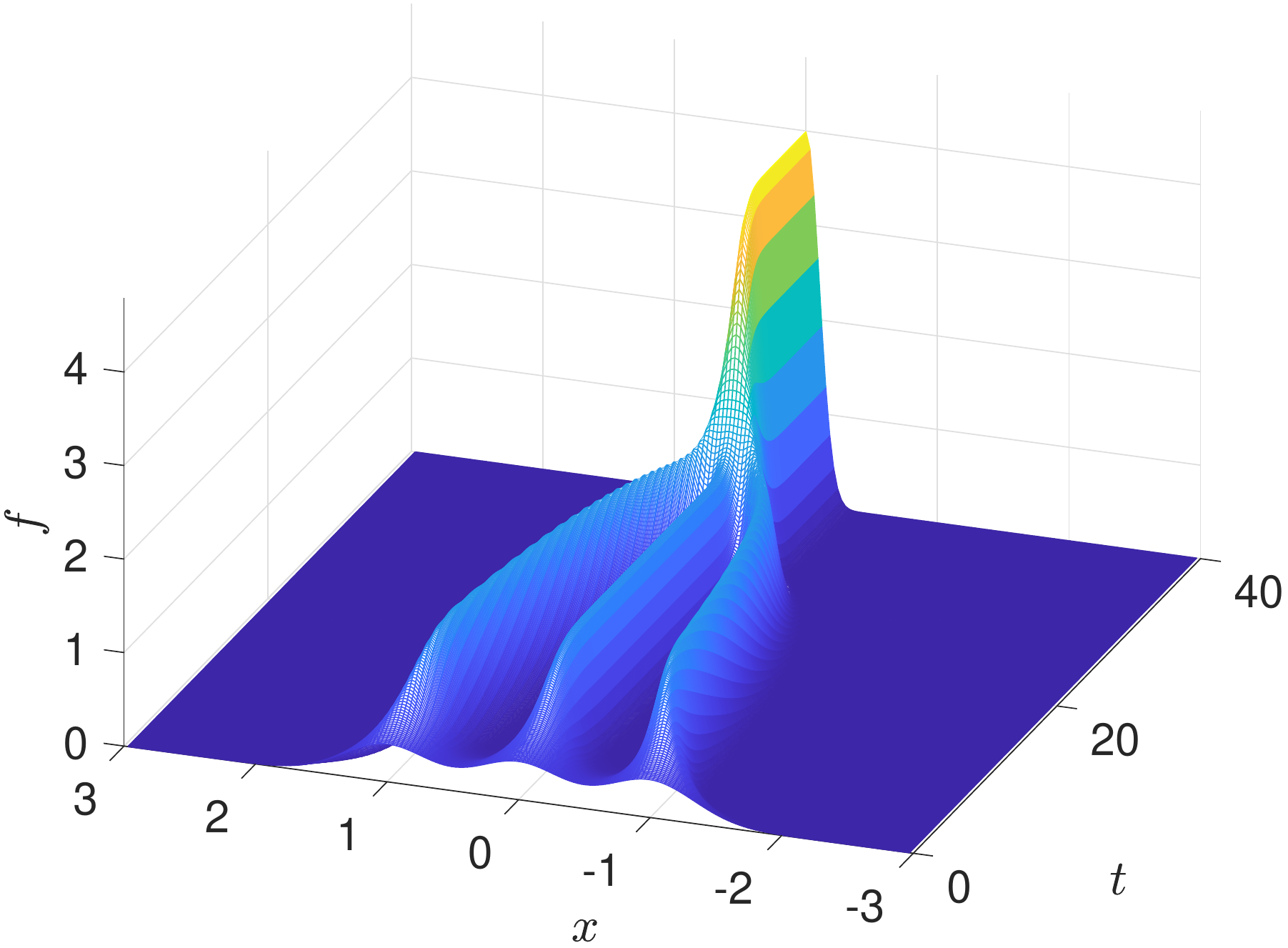}
\caption{Tightening and eventual collapse of three initial clusters.}
\label{fig:THREE_BUMP_3D_PLOTS}
\ec
\end{figure}

For the initial opinion distribution
$$
f_0(x) = \frac{1}{3} \left( \frac{e^{-5(x+1)^2}}{\sqrt{ \pi/5}} + \frac{e^{-5x^2}}{\sqrt{ \pi/5}} + \frac{e^{-5(x-1)^2}}{\sqrt{ \pi/5}} \right)
$$
we obtain the time evolution shown in Fig.\ \ref{fig:THREE_BUMP_3D_PLOTS}. 
A feature of some interest is that
the three clusters start out with equal amplitude, but by time 10, the middle cluster has a lower amplitude than the outlying
ones. A closer inspection of the computed density shows that 
this effect is mostly attributable  to less tightening in the central cluster, not to migration of individuals out of the central cluster: 
The percentage of individuals between $x=-0.5$ 
and $x=0.5$ is nearly exactly the same at time 10 as at time $0$. At approximately time 30, the three clusters merge into one.
This calculation was carried out with a bit less resolution: $m=100$ (so 199 Dirac delta functions), again $\Delta x=3/m$, and $\Delta t=0.1$, and $\sigma=0.1$. 

 \subsection{Concentration of the opinion holder density.} 
We saw in our numerical results that  multiple clusters in our model always appear
to be transient, eventually collapsing into a single
cluster. We will give the easy proof that this must always happen.

\vskip 10pt
\begin{proposition}
Assume that the $X_i$, $1 \leq i \leq n$, obey eq.\ (\ref{time_evolution}). Then 
$\min_{1 \leq i \leq n} X_i$ is increasing, $\max_{1 \leq i \leq n} X_i$ is decreasing, and 
$$
\lim_{t \rightarrow \infty}  \left(  \max_{1 \leq i \leq n} X_i - \min_{1 \leq i \leq n} X_i \right) = 0.
$$
\end{proposition}

\begin{proof} Without loss of generality, assume $X_1 = \min_{1 \leq i \leq n} X_i$ and $X_n = \max_{1 \leq i \leq n} X_i$.
(This is merely a matter of notation.)
We have 
\begin{eqnarray}
\nonumber
\frac{dX_1}{dt}  &=& \sum_{j=1}^n ~\!  \frac{ \eta( \left| X_j-X_1 \right|) w_j }{\sum_{\ell=1}^n {\eta( \left| X_\ell-X_1 \right|)w_\ell}} ~\! (X_j-X_1) \\
\nonumber
&\geq&  \frac{ \eta( \left| X_n-X_1 \right|) w_n }{\sum_{\ell=1}^n {\eta( \left| X_\ell-X_1 \right|)w_\ell}} ~\! (X_n-X_1)  \\
\nonumber
&\geq& \frac{ \eta( \left| X_n-X_1 \right|) w_n }{\sum_{\ell=1}^n {\eta( \left| X_n-X_1 \right|)w_\ell}} ~\! (X_n-X_1) \\
\label{eq:X_1_rises}
&=& w_n (X_n-X_1).
\end{eqnarray}
Similarly, 
\begin{eqnarray}
\nonumber
\frac{dX_n}{dt}  &=& \sum_{j=1}^n ~\!  \frac{ \eta( \left| X_j-X_n \right|) w_j }{\sum_{\ell=1}^n {\eta( \left| X_\ell-X_n \right|)w_\ell}} ~\! (X_j-X_n) \\
\nonumber
&\leq&  \frac{ \eta( \left| X_1-X_n \right|) w_1 }{\sum_{\ell=1}^n {\eta( \left| X_\ell-X_1 \right|)w_\ell}} ~\! (X_1-X_n)  \\
\nonumber
&\leq& \frac{ \eta( \left| X_1-X_n \right|) w_1 }{\sum_{\ell=1}^n {\eta( \left| X_1-X_n \right|)w_\ell}} ~\! (X_1-X_n) \\
\label{eq:X_n_falls}
&=& -w_1 (X_n-X_1).
\end{eqnarray}
Subtracting (\ref{eq:X_1_rises}) from (\ref{eq:X_n_falls}), we find
$$
\frac{d}{dt} (X_n-X_1) \leq -(w_1+w_n) (X_n-X_1).
$$
Since all the weights $w_j$  are positive, this implies the assertion.
\end{proof}
\vskip 10pt
For much less straightfoward results about convergence to consensus in {\em higher-dimensional}
 spaces, see 
\cite[Section 2]{Montsch_Tadmor}.

\section{A space-time-continuous model}
\subsection{Differential-integral formulation.} 
Let now $f=f(x,t)$ be a continuous opinion holder density. The analogues of 
$$
\sum_{j=1}^n \eta \left(  |X_i-X_j| \right) w_j X_j ~~~\mbox{and} ~~~ 
\sum_{\ell=1}^n \eta \left( | X_i-X_\ell|  \right) w_\ell
$$
are
$$
\int_{-\infty}^\infty \eta(|z|) f(x-z,t) (x-z) ~\! dz  ~~~\mbox{and} ~~~  \int_{-\infty}^\infty \eta(|z|) f(x-z,t) ~\! dz .
$$

We will derive a continuous evolution equation using arguments similar to those often used to 
derive conservation equations such as convection or diffusion equations.
Consider an interval $[a,b]$. At $a$, opinion holders are moving right with velocity 
$$
\frac{\int_{-\infty}^\infty \eta(|z|) f(a-z,t) (a-z) ~\! dz }{\int_{-\infty}^\infty \eta(|z|) f(a-z,t) ~\! dz} - a ~~= ~~
- \frac{\int_{-\infty}^\infty \eta(|z|) f(a-z,t) z ~\! dz }{\int_{-\infty}^\infty \eta(|z|) f(a-z,t) ~\! dz} .
$$
At $b$, they are similarly moving right with velocity 
$$
-\frac{\int_{-\infty}^\infty \eta(|z|) f(b-z,t) z ~\! dz }{\int_{-\infty}^\infty \eta(|z|) f(b-z,t) ~\! dz}.
$$
Now think about a short time interval of duration $\Delta t$. The fraction of opinion holders entering $[a,b]$ through $a$
in the time interval $[t,t+\Delta t]$ is about
$$
- f(a,t)  ~\! \frac{\int_{-\infty}^\infty \eta(|z|) f(a-z,t) z ~\! dz }{\int_{-\infty}^\infty \eta(|z|) f(a-z,t) ~\! dz}  ~ \Delta t.
$$
The fraction exiting through $b$ is similarly 
$$
- f(b,t)  ~\! \frac{\int_{-\infty}^\infty \eta(|z|) f(b-z,t) z ~\! dz }{\int_{-\infty}^\infty \eta(|z|) f(b-z,t) ~\! dz}  ~ \Delta t.
$$
It follows that 
$$
\frac{d}{dt} \int_a^b f(x,t) ~\! dx= 
\int_a^b f_t(x,t) ~\! dx  = 
 f(b,t)  ~\! \frac{\int_{-\infty}^\infty \eta(|z|) f(b-z,t) z ~\! dz }{\int_{-\infty}^\infty \eta(|z|) f(b-z,t) ~\! dz}  
$$
$$
\hskip 145pt  -~\! f(a,t)  ~\! \frac{\int_{-\infty}^\infty \eta(|z|) f(a-z,t) z ~\! dz }{\int_{-\infty}^\infty \eta(|z|) f(a-z,t) ~\! dz} .
$$
Since this holds for any choice of $[a,b]$, we conclude:

\begin{equation}
\label{differential_integral_equation}
f_t(x,t)  = \left( \frac{\int_{-\infty}^\infty \eta(|z|) f(x-z,t)  ~\!\! z~\!  dz }{\int_{-\infty}^\infty \eta(|z|) f(x-z,t) ~\! dz}  ~
f(x,t) \right)_x.
\end{equation}
\vskip 5pt
\noindent
The particle model presented in Section \ref{sec:particle_model} can
 be viewed as a numerical method for solving  eq.\ (\ref{differential_integral_equation}), with initial condition $f(x,0)=f_0(x)$ and 
 zero boundary conditions at $\pm \infty$.
 
 Equation (\ref{differential_integral_equation}) is closely related to others that have
 appeared in the literature, for instance \cite[eqs.\ (7) and (8)]{Canuto_2012}. 
We take the velocity at opinion
space location $x$ at time $t$ to be 
  \begin{equation}
 \label{eq:ours}
-  \frac{\int_{-\infty}^\infty \eta(|z|)  f(x-z,t) z~\! dz}{\int_{-\infty}^\infty \eta(|z|)  f(x-z,t) ~\! dz}  = \frac{\int_{-\infty}^\infty \eta(|y-x|) f(y,t)(y-x)  ~\! dy}{\int_{-\infty}^\infty \eta(|y-x|)  f(y,t)~\!  dy}
\end{equation}
while the velocity in \cite[eq.\ (8)]{Canuto_2012}, in the same notation, is
\begin{equation}
\label{eq:theirs}
\int_{-\infty}^\infty \eta(|y-x|)  f(y,t) (y-x) ~\! dy.
\end{equation}
In (\ref{eq:ours}), a weighted average of the differences $y-x$ is taken, while (\ref{eq:theirs}) will be larger if $x$ is 
surrounded by many nearby agents, smaller if it isn't. Which is more accurate depends on how 
opinion dynamics work --- we assume that all agents are equally eager to conform, even those surrounded by 
only few other agents, whereas \cite[eq.\ (8)]{Canuto_2012} implicitly assumes that those surrounded by many
agents are more eager to conform than those surrounded by few agents.

\subsection{Differential formulation.} 
In our numerical experiments, we always use
$
\eta(z) = e^{-z/\nu}
$
with $\nu>0$. This assumption was not crucial
until now, but will be here; we could not do the following computation with 
a general $\eta$.
There are two integrals in eq.\ (\ref{differential_integral_equation}), 
\begin{equation}
\label{eq:defg}
g(x,t) = \int_{-\infty}^\infty e^{-|z|/\nu} f(x-z,t) ~\! dz
~~~\mbox{and} ~~~ h(x,t)  = \int_{-\infty}^\infty z e^{-|z|/\nu} f(x-z,t) ~\! dz.
\end{equation}
The strategy is to write eq.\ (\ref{differential_integral_equation}) as
$$
f_t(x,t) = \left( \frac{h(x,t)}{g(x,t)}  ~ f(x,t) \right)_x
$$
and 
then add supplementary differential equations for $g$ and $h$.
We note  for later reference that 
\begin{equation}
\label{eq:total_mass_of_g}
\int_{-\infty}^\infty g(x,t) ~\! dx = 2 \nu  ~~~\mbox{for all $t$}.
\end{equation}
In fact, 
$
g(\cdot, t)  =  \eta \ast f(\cdot, t)
$
and $ \frac{\eta}{2 \nu}$ is a probability density. We refer to $g$ as the {\em locally averaged opinion holder density}. On the other hand, 
\begin{equation}
\label{eq:total_mass_of_h}
\int_{-\infty}^\infty h(x,t) ~\! dx= 0 ~~~\mbox{for all $t$}.
\end{equation}
The ratio 
$h(x,t)/g(x,t)
$
is a weighted average of positions, with a weight that is larger for views that are more commonly held and for views that
are closer to $x$.

\subsubsection{Supplementary differential equation for $\mathbf{g}$.} 
 We differentiate $g$ with respect to $x$, assuming
sufficient smoothness of $f$, and
use  
\begin{equation}
    \label{eq:xy}
\frac{\partial}{\partial x} f(x-z,t) = - \frac{\partial}{\partial z}  f(x-z,t) .
\end{equation}
We obtain: 
$$
g_x(x,t) = \int_{-\infty}^\infty e^{-|z|/\nu}  \frac{\partial}{\partial x} f(x-z,t) ~\! dz = 
-  \int_{-\infty}^\infty e^{-|z|/\nu}  \frac{\partial}{\partial z} f(x-z,t) ~\! dz.
$$
We integrate by parts to move the $z$-derivative to the exponential term, using the following formula, which
will be used several times in this section:
\begin{equation}
\label{eq:diff_exp}
\frac{d}{dz} \left( e^{-|z|/\nu} \right) = - \frac{1}{\nu}~\! \sign(z) ~\! e^{-|z|/\nu}.
\end{equation}
We find
\begin{equation}
\label{eq:g_x}
g_x(x,t) = - \frac{1}{\nu} \int_{-\infty}^\infty \sign(z) e^{-|z|/\nu} f(x-z,t) ~\! dz.
\end{equation}
We differentiate with respect to $x$ again, using this formula: 
\begin{equation}
\label{eq:diff_sign_exp}
\frac{d}{dz}  \left( \sign(z)  e^{-|z|/\nu}   \right) =2 \delta(z) - \frac{e^{-|z|/\nu}}{\nu}
\end{equation}
where $\delta$ is the Dirac delta function. This follows from the product rule, which is in fact rigorously applicable here, and
from the fact that $\sign(z) \cdot \sign(z) = 1$ for all $z \neq 0$.
We find
$$
g_{xx}(x,t) = -\frac{2}{\nu} f(x,t) + \frac{1}{\nu^2} g(x,t).
$$
We re-arrange a bit like this:
$$
- g_{xx} + \frac{g}{\nu^2}  =  \frac{2}{\nu} f.
$$
This is the supplementary differential equation for $g$.

\subsubsection{Supplementary differential equation for $\mathbf{h}$.} 
We will use similar reasoning for $h$, and for this purpose we first record that 
\begin{equation}
\label{eq:diff_z_exp} 
\frac{d}{dz} \left( z e^{-|z|/\nu} \right) = e^{-|z|/\nu} - \frac{1}{\nu} |z| e^{-|z|/\nu}.
\end{equation}
This follows from the product rule, using also that $z~\! \sign(z) = |z|$. Differentiation under the integral sign, using (\ref{eq:xy}), 
then integration by parts, now yields
$$
h_x(x,t) = \int_{-\infty}^\infty \left( e^{-|z|/\nu} - \frac{1}{\nu} |z| e^{-|z|/\nu} \right) f(x-z,t) ~\! dz = 
$$
$$
g(x,t) - \frac{1}{\nu} 
\int_{-\infty}^\infty |z| e^{-|z|/\nu} f(x-z,t) ~\! dz.
$$
We differentiate again with respect to $x$, using the following
 twin of formula (\ref{eq:diff_z_exp}):  
\begin{equation}
\label{eq:diff_abs_z_exp} 
\frac{d}{dz} \left( |z| e^{-|z|/\nu} \right) = \sign(z) e^{-|z|/\nu} - \frac{1}{\nu} z e^{-|z|/\nu}.
\end{equation}
We obtain, also using eq.\ (\ref{eq:g_x}): 
$$
h_{xx}(x,t) = g_x(x,t) - \frac{1}{\nu} \int_{-\infty}^\infty \sign(z) e^{-|z|/\nu} f(x-z) d z + \frac{1}{\nu^2} h(x,t) = 2 g_x(x,t) + 
\frac{1}{\nu^2} h(x,t).
$$
We re-arrange a bit to obtain
$$
-h_{xx} + \frac{1}{\nu^2} h = - 2 g_x.
$$
This is the supplementary differential equation for $h$.

\subsubsection{Summary of the differential formulation.} 
We have re-written
 eq.\ (\ref{differential_integral_equation}) as follows:
\begin{eqnarray}
\label{eq1}
f_t &=& \left( \frac{h}{g}  ~ f \right)_x, \\
\label{eq2}
- g_{xx} + \frac{g}{\nu^2}  &=&  \frac{2}{\nu} f, \\
\label{eq3}
-h_{xx} + \frac{h}{\nu^2}  &=& - 2 g_x.
\end{eqnarray}

\subsubsection{Concentration of the locally averaged opinion holder density.} 
We define
$$
H(x,t) = \int_{-\infty}^x h(s,t) ~\! ds.
$$
Note $H_x=h$ and $H(-\infty, t)=0$.
Integrating eq.\ (\ref{eq3}), we obtain
$$
- H_{xx} + \frac{H}{\nu^2} = -2 g + C
$$
with $C$ independent of $x$. 
Assuming  that $H_{xx} = h_x$ and $g$ vanish at $x=-\infty$, we conclude $C=0$, so 
\begin{equation}
\label{eq4}
-H_{xx} + \frac{H}{\nu^2} = - 2 g.
\end{equation}

Multiply both sides of (\ref{eq1}) by $H$, integrate with respect to $x$, and then integrate by parts on the right-hand side: 
\begin{equation}
\label{eq:rhs_done}
\int_{-\infty}^\infty H(x,t) f_t(x,t) ~\! dx = - \int_{-\infty}^\infty \frac{h^2(x,t)}{g(x,t)} f(x,t) ~\! dx.
\end{equation}
We'll re-write the
left-hand side of eq.\ (\ref{eq:rhs_done}) now. First, use eq.\ (\ref{eq2}): 
$$
\int_{-\infty}^{\infty} H(x,t) f_t(x,t) ~\! dx = \int_{-\infty}^\infty H(x,t)  \left( - \frac{\nu}{2} g_{xx} + \frac{g}{2 \nu} \right)_t ~\! dx.
$$
Integrating by parts twice, we obtain
$$
\int_{-\infty}^\infty \left( - \frac{\nu}{2} H_{xx} + \frac{1}{2 \nu} H \right) g_t ~\! dx, 
$$
and using (\ref{eq4}), this is 
$$
\int_{-\infty}^\infty - \nu g(x,t)  g_t(x,t) ~\! dx = - \frac{\nu}{2}  \frac{d}{dt} \int_{-\infty}^\infty g(x,t)^2 ~\! dx.
$$
In summary, canceling minus signs, (\ref{eq:rhs_done}) becomes
\begin{equation}
  \frac{d}{dt} \left\| g \right\|^2_{L^2} = \frac{2}{\nu} ~\!  \int_{-\infty}^\infty \frac{h^2(x,t)}{g(x,t)} f(x,t) ~\! dx.
 \end{equation}
 Positivity of $f$ implies positivity of $g$, and therefore 
 \begin{equation}
 \label{eq:increasing_L_2_norm} 
 \frac{d}{dt} \| g \|_{L^2}^2 > 0.
 \end{equation}
 Recall from eq.\ (\ref{eq:total_mass_of_g}) that the $L^1$-norm of $g$ is equal to $2 \nu$ for all time.  The square of the
 $L^2$-norm is a measure of concentration of $g$. 
This is reflected by the fact that if $X$ is a random number with probability density $\frac{g}{2 \nu}$, then 
 $$
 \| g \|_{L^2}^2 =  2 \nu E( g(X)),
 $$ 
 and therefore $E(g(X))$ rises as $t$ increases. When $X$ is drawn with density $g$, the expected value of $g(X)$ gets larger 
 with time. 
This means that $g$ becomes increasingly concentrated.

\section{Numerical convergence tests.} 
As an example, we test convergence for the initial condition
$$
f_0(x) = \frac{1}{3} \left( \frac{e^{-5(x+1)^2}}{\sqrt{ \pi/5}} + \frac{e^{-5x^2}}{\sqrt{ \pi/5}} + \frac{e^{-5(x-1)^2}}{\sqrt{ \pi/5}} \right).
$$
We track the $X_j$ up to time $t=1$. 
We use a time step $\Delta t$ and assume that $1/\Delta t$ is an integer. We initialize the $X_j$ at 
$$
X_j(0) = - 3 + j \Delta x, ~~~ j=1,2,\ldots,\frac{6}{\Delta x}-1,
$$
assuming that $6/\Delta x$ is an integer. 

\subsection{Weak convergence of the  computed opinion holder density.}
\label{subsec:weak_convergence_as_dx_and_dt_go_to_zero}
We  compute approximations for $f(x,t)$, $x = j \Delta x$, $j$ integer, $-3 < x < 3$, using
eq.\ (\ref{smooth}), where $\sigma = 0.1$. We denote these approximations by $f_{\Delta x, \Delta t}(x,t)$, and will 
test whether they converge to {\em some} limit as $\Delta x$ and $\Delta t$ are simultaneously reduced.
We cannot test convergence to an exact solution, since we have no 
analytic expression for an exact solution.

Fixing $\sigma$ independently of $\Delta x$  amounts to testing for a form
of  {\em weak} convergence. The computed distribution is 
a sum of $\delta$-functions, but we test for convergence of  the convolution with a Gaussian.

We define
$$
E_{\Delta x, \Delta t} = \max \left\{  \left| f_{\frac{\Delta x}{2}, \frac{ \Delta t}{2}}(x,1) - f_{\Delta x, \Delta t}(x,1) \right| ~:~ x = j \Delta x, ~ j ~\! \mbox{integer}, ~ -3 < x < 3 \right\}.
$$
If there is second-order convergence as $\Delta x$ and $\Delta t$ simultaneously tend to zero, one should expect 
$$
\frac{E_{\Delta x, \Delta t}}{E_{\Delta x/2, \Delta t/2}} \approx 4
$$
for small $\Delta x$ and $\Delta t$.
Table \ref{table:CONVERGENCE} 
confirms that this is indeed the case.

\begin{table}[h!]
\centering
\begin{tabular}{|c||c|c|c|}
\hline
$\Delta x$ & 0.06 & 0.03 &  0.015 \\
\hline 
$\Delta t $& 0.1 & 0.05 & 0.025 \\
\hline
$E_{\Delta x, \Delta t}$/$E_{\Delta x/2, \Delta t/2}$
&  4.01 & 3.98 & 4.00 \\
\hline
\end{tabular}
\caption{Numerical test confirming second order convergence of the approximation obtained by convolving the 
computed sum of delta functions with a Gaussian, as both $\Delta x$ and $\Delta t$ are refined; see text for 
details.} 
\label{table:CONVERGENCE}
\end{table}

%\subsection{Convergence of the computed g} 
%
%In Section \ref{subsec:weak_convergence_as_dx_and_dt_go_to_zero}, we considered the convolution 
%of the computed 
%$
%\sum_{i=1}^n  w_i \delta(x-X_i)
%$
%with a fixed Gaussian with standard deviation $\sigma=0.1$, and examined how rapidly the convolution converges
%as $\Delta x \rightarrow 0$ and  $\Delta t \rightarrow 0$.
%Here we will convolve with $\eta(z) = e^{-|z|/\nu}$ instead, obtaining approximations for the function $g$ defined
%in eq.\ (\ref{eq:defg}), and examine how fast those approximations converge. 
%We define $E_{\Delta x, \Delta t}$ as before, except using the computed approximations for $g$ now. 
%The ratios are still approximately 4, although far less accurately than in
%Section \ref{subsec:weak_convergence_as_dx_and_dt_go_to_zero}.
%
%%\renewcommand{\arraystretch}{1.2}
%\begin{table}[h!]
%\centering
%\begin{tabular}{|c||c|c|c|}
%\hline
%$\Delta x$ & 0.06 & 0.03 &  0.015 \\
%\hline 
%$\Delta t $& 0.1 & 0.05 & 0.025 \\
%\hline
%$E_{\Delta x, \Delta t}$/$E_{\Delta x/2, \Delta t/2}$
%&  3.5530 & 4.5870 & 3.3746 \\
%\hline
%\end{tabular}
%\caption{Numerical test confirming second order convergence of the approximation obtained by convolving the 
%computed sum of delta functions with $\eta(x) = e^{-|z|/\nu}$, as both $\Delta x$ and $\Delta t$ are refined.}
%\label{table:CONVERGENCE_2}
%\end{table}

\subsection{Dependence of the error on the time step.}
The most accurate calculation underlying Table  \ref{table:CONVERGENCE} uses
$$
\Delta x = \frac{0.06}{2^4} = 0.00375, ~~~ \Delta t = \frac{0.1}{2^4} = 0.00625.
$$
To test the importance of $\Delta t$ for the overall accuracy, we compare the results of this 
computation with results obtained using the same $\Delta x$, but coarser $\Delta t$.
We define
$$
F_{\Delta t} = \max \left\{  \left| f_{\Delta x = 0.00375, \Delta t =0.00625}(x,1) - f_{\Delta x=0.00375, \Delta t}(x,1) \right| ~:~ x = 0.00375j, \right. ~~~~~~~~
$$
$$
\hskip 240pt
\left. ~ j ~\! \mbox{integer}, ~ -3 < x < 3 \right\}.
$$
Table \ref{table:CONVERGENCE_3} shows the dependence of $F_{\Delta t}$ on $\Delta t$, confirming second order 
convergence.

\renewcommand{\arraystretch}{1.2}

\begin{table}[h!]
\centering
\begin{tabular}{|c||c|c|c|c|}
\hline 
$\Delta t $& 0.1 & 0.05 & 0.025 & 0.0125 \\
\hline
$F_{\Delta t}$ 
& $2.04 \times 10^{-5}$  &$5.12 \times 10^{-6}$  & $1.23 \times 10^{-6}$ &$ 2.47 \times 10^{-7}$ \\
\hline
$F_{\Delta t}/F_{\Delta t/2}$ & 3.98 & 4.17 &  4.98  & \\
\hline
\end{tabular}
\caption{Error as a function of $\Delta t$; see text for details.}
\label{table:CONVERGENCE_3}
\end{table}

\subsection{Dependence of the error on the spatial mesh size.} To test the importance of $\Delta x$ for the overall accuracy, we compare the results of the fine
computation using $\Delta x = 0.00375$ and $\Delta t = 0.00625$ with results obtained using the same $\Delta t$, but coarser $\Delta x$.
We define
$$
G_{\Delta x} = \max \left\{  \left| f_{\Delta x = 0.00375, \Delta t =0.00625}(x,1) - f_{\Delta x, \Delta t = 0.00625}(x,1) \right| ~:~ x = j \Delta x, \right. ~~~~~~~~
$$
$$
\hskip 240pt
\left. ~ j ~\! \mbox{integer}, ~ -3 < x < 3 \right\}.
$$
Table \ref{table:CONVERGENCE_3} shows the dependence of $G_{\Delta x}$ on $\Delta x$, again indicating
second order convergence.

\renewcommand{\arraystretch}{1.2}

\begin{table}[h!]
\centering
\begin{tabular}{|c||c|c|c|c|}
\hline 
$\Delta x $& 0.06 & 0.03 & 0.015 & 0.0075  \\
\hline
$G_{\Delta x}$ 
& $9.58 \times 10^{-4}$ & $2.36 \times 10^{-4}$ &$5.65 \times 10^{-5}$  &$ 1.13 \times 10^{-5} $ \\
\hline
$G_{\Delta x}/G_{\Delta x/2}$ & 4.06 &4.18 &  5.00 & \\
\hline
\end{tabular}
\caption{Error as a function of $\Delta x$; see text for details.}
\label{table:CONVERGENCE_4}
\end{table}

\section{Summary and discussion.} We began with a time-continuous version of Hegselmann-Krause dynamics,
similar to equations that have been proposed in the literature previously, but with
weighted particles, which we think of as representing clusters of agents, not individuals. The weights have a numerical advantage --- instead of needing many agents in a part of opinion space populated by many opinion holders, we can use 
fewer but heavier particles. 

The time-continuous model suggests  a fully continuous macroscopic model, which we formulated first as a single integral-differential equation, then --- for the special
case of an exponential interaction function --- as a system of differential equations. 

The time-continuous model 
(discretized using the midpoint method) 
can be viewed as a particle method for the fully continuous model. 
We demonstrated, numerically, the 
second-order convergence  of this method in a weak sense, meaning that the convolution of the solution with a mollifier is computed
with second-order accuracy.

In our numerical computations, all opinion holders eventually arrive at consensus. 
This could be counter-acted by adding {\em diffusion} (spontaneous random small changes in opinions) 
in the model, as some  authors  have proposed (see for instance \cite{Ben_Naim_2005,Goddard_2022}). We have refrained from
doing that here because it would raise, in our context, the question how to incorporate diffusion in the particle method.
One possibility would be a method similar to Chorin's random walk method for
viscous fluid dynamics \cite{Chorin_1973,Goodman_1987}.

Our model starts with the original Hegselmann-Krause model \cite{hegsel_krause_2002}.
In a later paper \cite{Hegselmann_Krause_2005}, Hegselmann and Krause suggested that individuals might respond not to 
the arithmetic average (or, in our modification of the model, weighted arithmetic average) of opinions in their vicinity, but to 
a different kind of average --- geometric averages for instance. We have not yet thought about what would happen if
we followed this interesting suggestion in our model.

In future work, we plan to use the method presented in this paper to explore the interaction of candidate dynamics  with voter opinion dynamics. We have taken a first step in that direction in \cite{Borgers_et_al_candidate_dynamics}.

\small{

}

  \end{document}